\documentclass[10pt]{amsart}

\usepackage{amsthm,amssymb,amsmath,amsopn,amsfonts}
\usepackage{stmaryrd,calc,mathrsfs,enumerate}
\usepackage{tikz}\usetikzlibrary{matrix,arrows} 

\usepackage{hyperref} 
\hypersetup{
	colorlinks=true,
	citecolor=blue,
	linkcolor=blue
}

\newtheorem{thm}{Theorem}[section]
\newtheorem*{thm*}{Theorem}
\newtheorem{prop}[thm]{Proposition}
\newtheorem{lem}[thm]{Lemma}
\newtheorem{cor}[thm]{Corollary}

\theoremstyle{definition}
\newtheorem{definition}[thm]{Definition}
\newtheorem{example}[thm]{Example}

\theoremstyle{remark}
\newtheorem{remark}[thm]{Remark}

\numberwithin{equation}{section}

\newcommand{\bdef}{\begin{definition}}\newcommand{\ndef}{\end{definition}}
\newcommand{\bteo}{\begin{thm}}\newcommand{\nteo}{\end{thm}}
\newcommand{\bprop}{\begin{prop}}\newcommand{\nprop}{\end{prop}}
\newcommand{\brmk}{\begin{remark}}\newcommand{\nrmk}{\end{remark}}
\newcommand{\bcor}{\begin{cor}}\newcommand{\ncor}{\end{cor}}
\newcommand{\blem}{\begin{lem}}\newcommand{\nlem}{\end{lem}}
\newcommand{\bex}{\begin{example}}\newcommand{\nex}{\end{example}}
 \newcommand{\Z}{\mathbb{Z}}
 
\newcommand{\C}{\mathbb{C}}  

\renewcommand{\O}{\mathcal{O}}

\newcommand{\g}{\mathcal}

\DeclareMathOperator{\Pic}{Pic}
\DeclareMathOperator{\Div}{Div}

\DeclareMathOperator{\cliff}{Cliff}
\DeclareMathOperator{\gon}{gon}

\begin{document}


\title[]{Gonality and Clifford index of curves on elliptic K3 surfaces with Picard number two}

\author{Marco Ramponi}
\address{Marco Ramponi, Universit\'{e} de Poitiers, Poitiers, France}
\email{marco.ramponi@math.univ-poitiers.fr}


\begin{abstract}
We compute the Clifford index of all curves on K3 surfaces with Picard group isomorphic to $U(m)$.
\end{abstract}

\maketitle

\section{Introduction}
In the past years many authors have studied problems related to the gonality and Clifford index of curves lying on K3 surfaces. In this paper, by \emph{curve} we always mean a smooth, reduced and irreducible curve over the field of complex numbers. The gonality and the Clifford index of a curve $C$ are respectively defined by
\begin{align*}
\gon(C) &= \min\{ \deg(A) :  A\in\Div(C), \ h^0(A)=2\}, \\
\cliff(C) &= \min\{ \cliff(A) :  A\in\Div(C), \ h^0(A)\geq2, \  h^1(A)\geq2 \}
\end{align*} 
where $\cliff(A)=\deg A-2h^0(A)+2$. The Clifford index measures how special $C$ is in the moduli space $\g{M}_g$ of curves of genus $g$, in the following sense. One has
$$0\leq\cliff(C)\leq \lfloor \tfrac{g-1}{2} \rfloor,$$
where the second inequality is an equality for the generic member of $\g{M}_g$, and on the other hand $\cliff(C)=0$ if and only if $C$ is hyperelliptic (cf.\! \cite{ACGH}).
Therefore, in this paper we say that a curve $C$ is \emph{Clifford general} if $\cliff(C)=\lfloor\tfrac{g-1}{2}\rfloor$. In all other cases, we say that $C$ is \emph{Clifford special}.

A classical result by Saint-Donat \cite[5.8]{SD} states that given a hyperelliptic curve $C$ on a K3 surface, all the curves in the linear system $|C|$ are also hyperelliptic. This interesting fact was vastly generalized by Green and Lazarsfeld \cite{GL}, who proved that indeed the Clifford index is the same for each smooth member of $|C|$. They also proved that, whenever $C$ is Clifford special, there exists a divisor $D$ on the ambient K3 surface $X$ whose restriction to $C$ computes its Clifford index\footnote{A divisor $A\in\Div(C)$ is said to compute the Clifford index of $C$ if it appears in the definition of $\cliff(C)$ and achieves the minimal value, i.e.\! $\cliff(A)=\cliff(C)$.}. In this case one says that \emph{the Clifford index is cut out on $C$ by $D$}. 
In \cite{Kn01} Knutsen showed that, if $C$ is Clifford special, we can choose $D$ to be a (smooth and irreducible) curve, and moreover $\cliff(C)=D.(C-D)-2$.

Originally, the constancy of the gonality of curves in a linear system on a K3 surface had been conjectured (unpublished) by Harris and Mumford, but Donagi and Morrison \cite{DM} found a counterexample. However, by the works of Ciliberto and Pareschi \cite{CP95} and Knutsen \cite{Kn09} we know that this is indeed the only counterexample.

On the other hand, the notions of gonality and Clifford index are very much related: for any curve $C$ of Clifford index $c$ one has
$$c+2\leq\gon(C)\leq c+3,$$
and curves for which $\gon(C)=c+3$ are conjectured to be very rare (cf.\! \cite{ELMS}). When lying on K3 surfaces, these curves are completely classified by Knutsen \cite{Kn09}.

In general, the gonality and the Clifford index are subtle invariants which are hard to compute explicitly for a given curve. In this note, we compute the Clifford index and gonality of all curves on some elliptic K3 surfaces. We prove the following.

\bteo\label{main-thm} Let $X$ be a K3 surface with Picard group isomorphic to $U(m)$, with $m\in\Z$, $m\ge1$. Denote by $E$ and $F$ two generators of  $\Pic(X)$, with $E^2=F^2=0$ and $E.F=m$. Let $C$ be a curve on $X$ of genus $g>2$. Then, either 
\begin{enumerate}[(i)]
\item The Clifford index of $C$ is cut out on $C$ by an elliptic curve $E_C$, which is linearly equivalent to the one among $E$ and $F$ having minimal intersection with $C$. Then $\cliff(C)=C.E_C-2$ and $C.E_C$ is equal to the gonality of $C$; or
\item $m>2$ and $C$ is linearly equivalent to $E+F$. Then $C$ has maximal Clifford index $\cliff(C)=\lfloor m/2\rfloor$.
\end{enumerate}
\nteo

In the statement of the theorem $U$ denotes the \emph{hyperbolic lattice}: the lattice given by $\Z\oplus\Z$ with intersection matrix
$$U=\begin{pmatrix} 0 & 1 \\ 1 & 0\end{pmatrix}.$$
$U(m)$ denotes the lattice obtained by $U$ by multiplying the intersection matrix by a non-zero integer $m$.  
The term \emph{isomorphic} means isomorphic as lattices.

Notice that there always exists a class of square zero in $\Pic(X)\simeq U(m)$. Therefore a K3 surface $X$ as in the theorem admits an elliptic fibration by \cite[\S 3]{ShapShaf}.

Part of the motivation for studying this problem came from a recent paper by Watanabe \cite{Wa} in which he shows that for a K3 surface $X$ which is a double cover of a smooth del Pezzo surface of degree $4\leq d\leq 8$ such that $X$ carries a non-symplectic automorphism of order two which acts trivially on the Picard group, then for any curve $C$ on $X$, either the Clifford index is cut out on $C$ by some elliptic curve on $X$, or $C$ is linearly equivalent to a multiple of a curve of genus 2. 
The key idea of Watanabe is that the automorphism yields some useful geometric informations which help to characterize the topological properties of the curves on $X$. 

In this work, we started to investigate the analogue situation when $X$ carries a non-symplectic automorphism of order 3 which acts trivially on the Picard group. For $\rho(X)=2$, we know by the classification results of Artebani-Sarti \cite{AS} and Taki \cite{tak11}, that the Picard group of $X$ is isomorphic to either $U$ or $U(3)$. Therefore, our result applies to this case. However, not all K3 surfaces with Picard group $U(m)$ admit non-symplectic automorphisms: see \cite{AS}, and also Artebani-Sarti-Taki \cite{AST}.

\subsection*{Notation and conventions}
We work over the complex number field $\C$. 
By \emph{surface} we mean a smooth irreducible projective surface. A K3 surface is a regular surface with trivial canonical bundle. The Picard number of a surface $X$ is by definition the rank of the Picard group $\Pic(X)$ and is denoted by the letter $\rho$.
The symbol $\sim$ denotes linear equivalence between divisors and $|D|$ is used to denote the complete linear system associated to a divisor $D$.
A lattice is a free $\Z$-module $L$ of finite rank equipped with a non-degenerate symmetric integral bilinear form $L\times L\rightarrow\Z$, $(x,y)\mapsto x.y$. An isomorphism of lattices is a $\Z$-module isomorphism preserving the bilinear forms.

\section{Clifford special curves on a K3 surface}\label{Clifford-theory}
Let $X$ be a K3 surface. In this short section we recall some fundamental results which will be needed in the following and we also explain why the case $\rho(X)=1$ is not interesting for our purposes. Fix a curve $C$ of genus $g$ on $X$. Let
$$\g{A}(C) := \{D\in\Div(X) \ | \ h^0(\O_X(D))\geq2, \ h^0(\O_X(C-D))\geq2\}.$$
Notice that $C$ admits a decomposition $C\sim D+D'$ into two moving classes $D$ and $D'$ if and only if $\g{A}(C)\neq\emptyset$. When this happens, among such decompositions it is interesting to consider those with minimal intersection $D.D'$. Hence, one defines
$$\mu_C := \min\{D.(C-D)-2 \ | \ D\in\g{A}(C)\},$$
and denotes by $\g{A}^0(C)$ the divisors in $\g{A}(C)$ achieving this minimal value:
$$\g{A}^0(C) := \{D\in\g{A}(C) \ | \ D.(C-D)-2 = \mu_C\}\subset\g{A}(C).$$

Observe that $\mu_C\geq0$ since the curve $C$ (or any member of the complete linear system of a basepoint free and big line bundle on a K3 surface) is numerically $2$-connected (cf.\! \cite[(3.9.6)]{SD}).
By the results in \cite{GL} and \cite{Kn01} we have: (cf.\! \cite[p.11]{Kn04})
\begin{equation*}
\cliff(C)=\min\{\mu_C, \lfloor \tfrac{g-1}{2} \rfloor\}.
\end{equation*}

In other words, either $C$ is Clifford general, or $C$ is Clifford special and then $\g{A}(C)$ is non-empty, the Clifford index of $C$ is cut out by some divisor $D\in\g{A}^0(C)$ on $X$ and $\cliff(C)=\mu_C$. Then, by definition of Clifford index
$$h^0(\O_C(D))\geq2  \mbox{ \ and \ }  h^1(\O_C(D))=h^0(\omega_C\otimes\O_C(-D))\geq2.$$
In particular, the linear systems on the curve $C$ given by the line bundles $\O_C(D)$ and $\omega_C\otimes\O_C(-D)=\O_C(C-D)$ contain some non-trivial effective divisors, hence of positive degree. Therefore we get $\deg_C(C-D)=C.(C-D)>0$ and also $\deg_C(D)=C.D>0$. Altogether, this yields the following inequalities:
\begin{equation*}
0<C.D<C^2.
\end{equation*}
In particular, when $\Pic(X)=\Z[H]$, we see that this inequalities are impossible when $C\in|H|$, and so $C$ has general Clifford index in this case. On the other hand, for $C\sim kH$ with $k\geq2$, a direct computation shows that the Clifford index of $C$ is cut out on $C$ by a member of $|H|$. When $\rho(X)\geq2$, however, the situation is more interesting. In the next section we compute the Clifford index of any curve on a K3 surface with Picard group isomorphic to $U(m)$.

\section{Proof of the Theorem}
Let $X$ be a K3 surface with Picard group isomorphic to $U(m)$, with $m\geq1$. 

We let the Picard group of $X$ be generated by the classes of two effective divisors $E$ and $F$ such that $E^2=F^2=0$ and $E.F=m$. Up to the action of the Weyl group of $X$ we may assume that $E$ is an elliptic curve (cf.\! \cite[\S 3]{ShapShaf}). 

When $m=1$, we observe that the rational curve $\Gamma\sim F-E$ yields a section of the elliptic fibration given by $|E|$. Moreover, the linear system $|F|=|E+\Gamma|$ contains a rational curve as a base component and therefore $F$ cannot be represented by an irreducible curve (cf.\! Saint-Donat \cite[2.6 \& 2.7]{SD}).

On the other hand, when $m>1$, since $x^2\in2m\Z$ for $x\in U(m)$ we observe that there are no rational curves on $X$. Thus any effective divisor is nef and basepoint free (ibid.). In particular we may assume that $F$ is an elliptic curve. A simple computation shows that any elliptic curve on $X$ belongs to either $|E|$ or $|F|$.

For any effective divisor $C$ on $X$ let us define
\[\begin{array}{lll}
d_C &:=& \min\{E' \cdot C \ | \ E' \mbox{ is an elliptic curve on } X\}, \\
\g{E}^0(C) &:=& \{\mbox{elliptic curves } E_C \mbox{ such that } E_C \cdot C=d_C\}.
\end{array}\]

\blem\label{E+F} Let $C$ be a curve with $C^2>0$ and let $E_C$ be an elliptic curve in $\g{E}^0(C)$. If $(C-E_C)^2=0$, then $C$ belongs to the linear system $|E+F|$.\nlem
\begin{proof}
If $(C-E_C)^2=0$ then $C-E_C$ is linearly equivalent to a multiple of an elliptic curve $E'$, so that we can write $C = E_C+(C-E_C)\sim E_C+kE'$, some $k\geq1$. Since $C^2>0$ we see that $E'$ is not linearly equivalent to $E_C$. Since $E'.C=E_C.C$, we get $k=1$ and $C\sim E_C+E'\sim E+F$.
\end{proof}

\blem\label{Clifford-General} Let $m\geq2$ and let $C$ be a curve in the linear system $|E+F|$. Then
\begin{enumerate}[(i)]
\item If $m=2$ then $C$ is Clifford special.
\item If $m>2$ then $C$ is Clifford general.
\end{enumerate}
\nlem

\begin{proof}
Assume $D\in\g{A}^0(C)$ and let $D\sim aE+bF$, with $a,b\geq0$. Then by definition of $\g{A}(C)$ we may assume $C-D$ effective, so that $0\leq(C-D).E=m(1-b)$ and $0\leq(C-D).F=m(1-a)$. Hence $a,b\in\{0,1\}$. This shows that the only curves $D$ in $\g{A}^0(C)$ are the members of $|E|$ and $|F|$. Then $C$ is Clifford special whenever
$\mu_C=D.(C-D)-2=m-2<\lfloor C^2/4\rfloor=\lfloor m/2\rfloor$, that is for $m\leq2$.
\end{proof}
\brmk The case $m=1$ is not to be considered here since the linear system $|E+F|$ contains a rational curve $\Gamma\sim F-E$ as base component; hence there are no (irreducible) curves in $|E+F|$ in this case.
\nrmk

\blem\label{minusE} Let $C\subset X$ be an effective divisor with $C^2>0$. For any elliptic curve $E'$ on $X$ we have 
$$(C-E')^2\geq0, \quad h^0(C-E')\geq2.$$ 
Moreover, $|C-E'|$ is basepoint free for $m>1$. 
\nlem
\begin{proof}
Since $E$ and $F$ are the only effective reduced divisors with self-intersection zero, it is clear that in order to show the Lemma we may assume $E'\in|E|$, by the symmetry of the roles of $E$ and $F$ in $\Pic(X)$. Let $C\sim aE+bF$ for some positive integers $a$ and $b$.
Then clearly $(C-E)^2\geq0$ and also $C.E>0$. Thus $E.(C-E)>0$, which shows that $C-E$ is effective. It follows $h^0(C-E)\geq2$ by Riemann-Roch. Moreover, if $m>1$, then $|C-E|$ is basepoint free since in this case there are no rational curves on $X$ (cf.~\cite[\S 2.7]{SD}).
\end{proof}

\brmk Let $C$ be a curve on $X$. By the definition of $\g{A}(C)$ and Lemma \ref{minusE} above $\g{E}^0(C)\subset\g{A}(C)$. In particular $\mu_C\leq d_C-2.$
Moreover, 
$$\g{E}^0(C)\subset\g{A}^0(C) \iff \mu_C=d_C-2.$$ 
Indeed, let $E_C\in\g{E}^0(C)\subset\g{A}(C)$. If $\mu_C=d_C-2$ then $E_C$ computes $\mu_C$. Hence $E_C\in\g{A}^0(C)$. The other implication is obvious.
\nrmk


\begin{proof}[Proof of Theorem \ref{main-thm}]
The first (and longer) part of the proof is to show that 
$$\g{E}^0(C)\subset\g{A}^0(C).$$
Let $E_C\in\g{E}^0(C)$. By Lemma \ref{minusE}, $E_C\in\g{A}(C)$, so that $\g{A}^0(C)$ is not empty. 

If $\g{A}^0(C)$ contains some elliptic curve $F$, then $C.E_C\leq C.F$ and so $C.E_C-2\leq\mu_C$. Since $E_C\in\g{A}(C)$, we have $C.E_C-2=\mu_C$. Therefore $E_C\in\g{A}^0(C)$. 

So we assume that $\g{A}^0(C)$ contains no elliptic curves at all. Let $D$ be an effective divisor in $\g{A}^0(C)$. Since $D\in\g{A}(C)$, and $h^1(D)=0$ by \cite[Prop.\! 2.6]{Kn04}, we have $D^2\geq0$. Let us show that, in fact, 
$$D^2\geq2.$$ 
Indeed, assume by contradiction $D^2=0$.
\begin{itemize}
\item In the case where $m\geq2$, since $X$ contains no rational curves, $D$ is basepoint free, and so it is linearly equivalent to an elliptic curve, by \cite[2.6]{SD}. This contradicts the assumption that $\g{A}^0(C)$ contains no elliptic curves.
\item In the case where $m=1$, let $E$ and $F$ be generators of the Picard group of $X$, with $E^2=F^2=0$ and $E.F=1$. Then we may assume that $E$ is an elliptic curve and there exists a rational curve $\Gamma$ on $X$ such that $\Gamma\sim F-E$. Since $D^2=0$ and $D\in\g{A}^0(C)$, by \cite[Prop.\! 2.6]{Kn04} we have $D\sim E$ or $F$. However, $E$ is not in $\g{A}^0(C)$ by assumption, thus $D\sim F$. Since $\Gamma$ is the base locus of $F$, we have $C.(F-E)=0$. On the other hand,
$$\mu_C=C \cdot D-2=C \cdot E-2$$
and, since $E\in\g{A}(C)$, this yields $E\in\g{A}^0(C)$. A contradiction.
\end{itemize}

By the above discussion we always have $D^2\geq2$. Notice that $C-D\in\g{A}^0(C)$ and then $(C-D)^2\geq2$ by the same reason. We want to show that $E_C\in\g{A}^0(C)$, contradicting the assumption that $\g{A}^0(C)$ contains no elliptic curves. Concretely, we need to show the following inequality
$$E_C \cdot C\leq D \cdot (C-D).$$ 
Rewrite this inequality as
\begin{equation}\label{inequality1}
(D-E_C)\cdot(D'-E_C)\geq0, \quad D':=C-D
\end{equation}
For $E_D\in\g{E}^0(D)$ and $E_{D'}\in\g{E}^0(D')$ we let
\begin{align*}
n_D &= (D-E_D).(D'-E_{D'}) \\
r_D &= D.(E_{D'}-E_C) \\
r_{D'} &= D'.(E_{D}-E_C)
\end{align*}
so that we may now rewrite (\ref{inequality1}) as follows:
\begin{equation}\label{inequality2}
n_D+r_D+r_{D'}\geq E_D \cdot E_{D'}
\end{equation}

\textbf{Claim.} For \emph{any choice} of elliptic curves $E_D\in\g{E}^0(D)$ and $E_{D'}\in\g{E}^0(D')$, 
$$n_D=(D-E_D).(D'-E_{D'})\geq0.$$ 
Indeed, by Lemma \ref{minusE} the classes of $(D-E_D)$ and $(D'-E_{D'})$ have non-negative self-intersection and are effective, thus they lie in the closure of the positive cone and intersect non-negatively (cf.\! \cite[IV.7]{BHPV}). This proves our claim.

Now, consider the following inequalities:
\begin{align}\label{r_D}
\begin{split}
r_D & \geq r_D+C.(E_C-E_{D'}) = D'.(E_C-E_{D'})\geq0 \\
r_{D'} &\geq r_{D'}+C.(E_C-E_{D}) = D.(E_C-E_D)\geq0
\end{split}
\end{align}

If we assume that either $r_D>0$ or $r_{D'}>0$ then (\ref{inequality2}) holds, since $n_D\geq0$ and
$$E_C \cdot E_{D'}\leq m \quad \mbox{and} \quad r_D\geq m \mbox{ \ or \ } r_{D'}\geq m.$$ 
(recall that $x.y\in m\Z$ for $x,y\in U(m)$). Hence, we assume $r_D=0$ and $r_{D'}=0$. Substituting this in (\ref{r_D}) we get $D.E_C=D.E_D$ and $D'.E_C=D'E_{D'}$. Thus, 
$$E_C\in\g{E}^0(D)\cap\g{E}^0(D').$$ 
Using the claim above, we can replace both $E_D$ and $E_{D'}$ by $E_C$ in the definition of $n_D$ and this yields the desired inequality (\ref{inequality1}) and therefore $\g{E}^0(C)\subset\g{A}^0(C)$.

Now that we know $\g{E}^0(C)\subset\g{A}^0(C)$, we determine all Clifford general curves. 
Take $E_C\in\g{E}^0(C)$. By Lemma \ref{minusE} we know $C-E_C\in\g{A}(C)$. Moreover, we also have $C-E_C\in\g{A}^0(C)$ since $E_C\in\g{A}^0(C)$ by assumption. We distinguish two cases:
\begin{itemize}
\item $(C-E_C)^2=0$. Then, by Lemma \ref{E+F} and Lemma \ref{Clifford-General}, $C$ is Clifford general if and only if $m>2$ and $C\in|E+F|$.
\item $(C-E_C)^2>0$. (in particular $C$ is not linearly equivalent to $E+F$). We then show that $C$ is Clifford special. This amounts to show 
$$2\mu_C\leq g-3$$ 
which, by the definition of $\mu_C$ and the genus formula, is equivalent to 
$$(C-2E_C)^2\geq -4.$$ 
We may write $C\sim aE_C+D$, with $a\geq1$, $D$ effective and $D^2=0$. If $a=1$ we get $C\in|E+F|$ by Lemma \ref{E+F}, which is not the case. So $a\geq2$ and 
$$(C-2E_C)^2\geq0.$$ 
Therefore, $C$ is Clifford special.
\end{itemize}

This proves that $C$ is Clifford general if and only if $m>2$ and $C\in|E+F|$, as in part (ii) of the Theorem. To show part (i), we can therefore assume that $C$ is Clifford special. Then $\cliff(C)=\mu_C$ and since $\g{E}^0(C)\subset\g{A}^0(C)$ we have $\mu_C=d_C-2$. Therefore, the Clifford index of $C$ is cut out by some elliptic curve $E_C\in\g{E}^0(C)$. In particular, $\cliff(C)$ is computed by a pencil: the restriction of $|E_C|$ to $C$. Therefore (cf.\! \cite[p.174]{ELMS}) 
$$\gon(C)=\cliff(C)+2=d_C.$$ 
Hence, the assertions of (i) follow and the Theorem is proved.
\end{proof}

\brmk In particular, we observe that when $m=1$ or $2$, any curve on $X$ is Clifford special and its Clifford index is cut out by an elliptic curve. The same conclusion when $m=2$ is implicitly contained in \cite{Wa}.
\nrmk

\subsection*{Acknowledgement}
I am thankful to Alessandra Sarti for her guidance and support. A special thank also to Flaminio Flamini and Andreas Knutsen for their warm welcome in Rome and useful conversations. I am grateful to Kenta Watanabe for a careful review and for pointing out some mistakes in a draft version of this paper.

\bibliographystyle{plain}

\bibliography{bigbib}

\end{document}